\documentclass[12pt,twoside]{article}
\usepackage{graphicx}
\usepackage{amsfonts,amsbsy,amssymb,amsmath, amsthm}
\allowdisplaybreaks
\usepackage{cite}
\usepackage{graphics}
\def\bpsp{\begin{pspicture}}
\def\epsp{\end{pspicture}}

\newtheorem{theorem}{Theorem}[section]
\newtheorem{remark}[theorem]{Remark}
\newtheorem{example}[theorem]{Example}
\newtheorem{lemma}[theorem]{Lemma}
\newtheorem{corollary}[theorem]{Corollary}
\newtheorem{definition}[theorem]{Definition}
\newtheorem{proposition}[theorem]{Proposition}

\newtheorem{note}{Note}
\newtheorem{case}{Case}

\newtheorem{conjecture}{Conjecture}
\newtheorem{question}{Question}

\newcommand{\bea}{\begin{eqnarray}}
\newcommand{\eea}{\end{eqnarray}}
\newcommand{\beq}{\begin{eqnarray*}}
\newcommand{\eeq}{\end{eqnarray*}}

\def\m4{\mbox{\rm ~(mod $4$)}}

\def \bd{\begin{definition}}
\def \ed{\end{definition}}
\def \bqu{\begin{question}}
\def \equ{\end{question}}
\def \bcc{\begin{conjecture}}
\def \ecc{\end{conjecture}}
\def \bt{\begin{theorem}}
\def \et{\end{theorem}}
\def \bl{\begin{lemma}}
\def \el{\end{lemma}}
\def \bc{\begin{corollary}}
\def \ec{\end{corollary}}
\def \be{\begin{equation}}
\def \ee{\end{equation}}
\def \ben{\begin{enumerate}}
\def \een{\end{enumerate}}
\def \ba{\begin{array}}
\def \ea{\end{array}}
\def \bp{\begin{proposition}}
\def \ep{\end{proposition}}
\def \bx{\begin{example}}
\def \ex{\end{example}}
\def \br{\begin{remark}}
\def \er{\end{remark}}
\def \bdsc{\begin{description}}
\def \edsc{\end{description}}

\def \bn{\begin{case}}
\def \en{\end{case}}
\def \bnt{\begin{note}}
\def \ent{\end{note}}
\def\1{1\!\!1}

\def\mm2{\mbox{\rm ~(mod $2$)}}
\def\m4{\mbox{\rm ~(mod $4$)}}

\def\m{\mu}

\def\1{\textbf{1}}
\def\0{\textbf{0}}

\begin{document}
\title{On the spectral redundancy of pineapple graphs }
\author{Pawan Kumar$^{a}$, S. Pirzada$^{b}$, Merajuddin$^{c}$, Yilun
Shang$^{d}$ \\
$^{a,c}${\em Department of Applied Mathematics, AMU, Aligarh, India}\\
$^{a,c}${\em Department of Mathematics, University of Kashmir, Srinagar, India}\\
$^{d}${\em Department of Computer and Information Sciences,}\\
{\em Northumbria University, Newcastle, UK}\\
$^{a}$\texttt{paone101@gmail.com}; $ ^{b} $\texttt{pirzadasd@kashmiruniversity.ac.in}\\
$^{c}$\texttt{meraj1957@rediffmail.com}; $ ^{d} $\texttt{yilun.shang@northumbria.ac.uk}\\
}

\date{}

\pagestyle{myheadings} \markboth{Kumar, Pirzada, Merajuddin, Shang }{On the spectral redundancy of pineapple graphs }
\maketitle
\vskip 5mm
\noindent{\footnotesize \bf Abstract.} In this article, we explore the concept of spectral redundancy within the class of pineapple graphs, denoted as $\mathcal{P}(\alpha,\beta)$. These graphs are constructed by attaching $\beta$ pendent edges to a single vertex of a complete graph $K_\alpha$. A connected graph $G$ earns the title of being spectrally non-redundant if the spectral radii of its connected induced subgraphs remain distinct. Spectral redundancy, on the other hand, arises when there is a repetition of spectral radii among the connected induced subgraphs within $G$.  Specifically, we analyze the adjacency spectrum of $\mathcal{P}(\alpha,\beta)$, revealing distinct eigenvalues including $0$, $-1$, and additional eigenvalues, some negative and others positive. Our investigation focuses on determining the spectral redundancy within this class of graphs, shedding light on their unique structural properties and implications for graph theory.

\vskip 3mm



\section{Introduction}

\noindent Suppose that $A$ is an $n$ by $n$ real matrix. A scalar
$\lambda \in R$ is a complementarity eigenvalue of $A$ if there is a
non-zero vector $ x \in {R}^n $ satisfying $ x \geq 0, \quad  Ax-\lambda x \geq 0$  and $\langle x,  Ax-\lambda x \rangle = 0$, where
$x\geq 0$ means $x$ is non-negative component wise. Seeger
\cite{intr}, introduced the complementarity spectrum of a square
matrix. Later, Fernandes et al. \cite{Fer} defined the
\emph{complementarity spectrum of a graph} $G$ to be the set of the
complementarity eigenvalues of the adjacency matrix $A_G$ of $G$. We
refer the readers to Seeger \cite{See2020,See2018}, Seeger and Sossa
\cite{Sos2021,see2,friendship}, Pinheiro et al. \cite{Pin2020} and Merajuddin et al. \cite{ppmt}
for more recent results on complementarity spectrum.

Let $G(V,E)$ be a simple connected graph, where $V= \{
v_1,v_2,\ldots , v_n\}$ is the vertex set and $E=\{e_1,e_2,\ldots
,e_m\}$ is the edge set. The graph $G$ has a symmetric adjacency
matrix $A(G)=(a_{ij})$ of order $n$, where

\[
    a_{ij}=
\begin{cases}
    1,& \text{if } v_i ~ \text{and } v_j ~\text{are~ adjacent}\\
    0,              & \text{otherwise.}
\end{cases}
\]
The largest root of the characteristic polynomial of the adjacency
matrix $A(G)$ is referred to as the the spectral radius of $G$.
Fernandes et al. \cite{Fer}, showed that the complementary spectrum
of $G$, denoted by $\Pi(G)$, is the set composed by the spectral
radius of the adjacency matrices of all the induced non-isomorphic
connected subgraphs of $G$, that is, $ \Pi(G)= \{\rho(F): F \in
\mathcal{S}(G) \}$. Here, $\mathcal{S}$($G$) signifies the
collection of all induced connected subgraphs of $G$. As a
consequence, if $b(G)$ and $c(G)$ are the number of all induced
non-isomorphic subgraphs and the cardinality of $\Pi(G)$,
respectively, the following inequalities are due to Fernandes et al.
\cite{Fer},
\begin{equation}\label{eq:bnd}
c(G)\leq b(G) \quad \text{and} \quad  n\leq b(G) \leq 2^n-1.
\end{equation}
The upper bound in the second inequality $b(G) \leq 2^n-1$ is not
tight. However, in \cite{Fer}, it is proven that it grows faster
than any polynomial in $n$. The lower bound $n\leq b(G)$ has been
nicely settled in \cite{See2018} and equality in this case holds
only for the \emph{elementary graphs:} paths, stars, cycles and
complete graphs. Any other graph with $n$ vertices has more than $n$
induced subgraphs. Since there may be (non-isomorphic) induced
subgraphs of a given graph $G$ with the same spectral radius, the
first inequality $c(G)\leq b(G)$ follows. In an attempt to measure
this repetition, Seeger~\cite{See2020} defined the \emph{spectral
redundancy} of a given graph $G$ as the ratio
$$ r(G) = \frac{b(G)}{c(G)}. $$
For a set $\mathcal{G}$ of connected graphs, the \emph{spectral
redundancy index} of $\mathcal{G}$ is defined as $r(\mathcal{G})=
\underset{G\in \mathcal{G}}{sup}\, \, r(G).$

The pineapple graph $\mathcal{P}(\alpha,\beta)$ is the coalescence
of the star $S_{1,\beta}$ at the vertex of degree $\beta$ with the
complete graph $K_\alpha$ at any vertex. The star graph
$S_{1,\beta+1}$ and the complete graph $K_\alpha$ are both special
cases of the pineapple graph when $\alpha=2$ and $\beta=0$,
respectively. If we include these two special cases in the family of
pineapple graphs, the family of pineapple graphs is hereditary, that
is, the induced subgraph of a pineapple graph is again a pineapple
graph. When $\beta>0$, all induced subgraphs with at least one
pendant vertex share the same unique groupie \cite{sha}.

In this article, we study the spectral aspect of pineapple graphs.
It was claimed in \cite{zhang2009some} that the pineapple graphs are
determined by their adjacency spectrum. But Topcu et al.
\cite{topcu2016spectral}, found some disconnected graphs with the
same characteristic polynomial as the pineapple graphs. However, it
is shown in \cite{topcu2019graphs} that the pineapple graph can
determined by its adjacency spectrum if we are confined in the realm
of connected graphs. The authors also determined all the
disconnected graphs which are cospectral with a pineapple graph.

The organization of the rest of the paper is as follows. In Section
2, we study the spectral radius of the pineapple graphs and give
some results regarding the spectral redundancy of this family.

\section{Spectral redundancy of pineapple graphs }

We begin with the following observation.

\begin{lemma}\label{coal}
The coalescence of two graphs $G$ and $H$ by identifying a vertex
$u$ of $G$ with a vertex $v$ of $H$ is denoted by $G\cdot H$. We have the
following relationship: $ P_{G\cdot H}(\lambda)=
P_{G}(\lambda)P_{H-v}(\lambda)+P_{G-u}(\lambda)P_{H}(\lambda)-\lambda
P_{G-u}(\lambda)P_{H-v}(\lambda).$

\end{lemma}

\begin{figure}
\centering
\includegraphics[width=5cm, height=6cm]{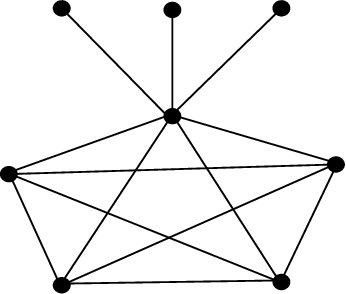}
\caption{Pineapple graph, $\mathcal{P}(5,3)$, with clique 5 and 3
pendant vertices}
\end{figure}

\noindent First, notice that the pineapple graph
$\mathcal{P}(2,\beta)$ is the star graph with $\beta+1$ pendant
vertices and the pineapple graph $\mathcal{P}(\alpha, 0)$ is the
complete graph isomorphic to $K_\alpha$. The following lemma is
about the spectral radius of the pineapple graph
$\mathcal{P}(\alpha,\beta)$ with $\alpha\geq 2$ and $\beta\geq 0$.

\begin{lemma}(\cite{topcu2016spectral})
    If $\mathcal{P}(\alpha,\beta)$ is the pineapple graph, then the characteristic equation of the adjacency matrix of $\mathcal{P}(\alpha,\beta)$ is given by
    $$P(\alpha,\beta,x)= x^{\beta-1}(x+1)^{\alpha-2}[x^3-(\alpha-2)x^2-(\alpha+\beta-1)x+\beta(\alpha-2)]=0.$$
\end{lemma}

\begin{figure}
\centering
\includegraphics[width=12cm, height=6.75cm]{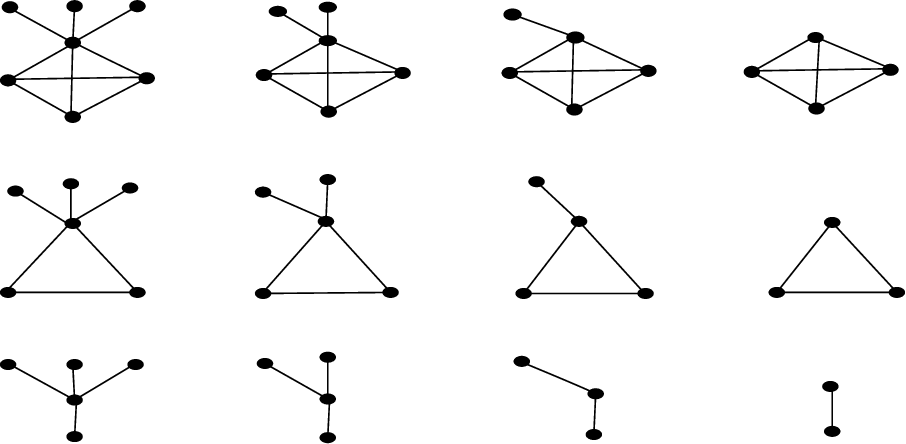}
\caption{ Non empty connected non-isomorphic induced subgraphs of
$\mathcal{P}(4,3)$}\label{indgraph}
\end{figure}

As the spectral radius of a connected graph is always positive, so
the spectral radius of the pineapple graph
$\mathcal{P}(\alpha,\beta)$ is the largest root of the equation
\begin{equation}\label{char}
   P(\alpha,\beta,x)=x^3-(\alpha-2)x^2-(\alpha+\beta-1)x+\beta(\alpha-2)=0
\end{equation}
We have the following useful lemma.
\begin{lemma}\label{lem23}
    For all $\alpha \geq 2$ and $\beta\geq 0$, we have $b(\mathcal{P}(\alpha,\beta))= (\alpha-1)(\beta+1)+1.$
\end{lemma}
\begin{proof} Clearly, an induced subgraph of $\mathcal{P}(\alpha,\beta)$ is again a pineapple graph (considering complete graph and star graph as pineapple graphs) with lesser number of vertices. There are exactly $\beta+1$ induced subgraph with fixed clique number $k$ $(2\leq k\leq \alpha)$. All of them are non-isomorphic. Therefore, the number of non empty connected induced subgraphs of $\mathcal{P}(\alpha,\beta)$ is $(\alpha-1)(\beta+1)$. There is an induced subgraph with one vertex as well. So $b(\mathcal{P}(\alpha,\beta))= (\alpha-1)(\beta+1)+1$. For the connected induced subgraph of $\mathcal{P}(4,3)$, see Fig.~$\ref{indgraph}$.
\end{proof}

The following lemmas will be useful.
\begin{lemma}\label{ind}
    If $H$ is a proper induced subgraph of $G$, then $\rho(H) < \rho(G)$
\end{lemma}
\begin{lemma}
     If $\mathcal{P}_1(\alpha_1,\beta_1)$ and $\mathcal{P}_2(\alpha_2,\beta_2)$ are two non-isomorphic pineapple graphs having the same spectral radius, then $(\alpha_2-\alpha_1)(\beta_2-\beta_1)<0$.
\end{lemma}
\begin{proof}
    Let $\mathcal{P}_1(\alpha_1,\beta_1)$ and $\mathcal{P}_2(\alpha_2,\beta_2)$ be two non-isomorphic pineapple graphs such that $\rho(\mathcal{P}_1(\alpha_1,\beta_1))= \rho(\mathcal{P}_2(\alpha_2,\beta_2))$. Without loss of generality, let us assume $(\alpha_2-\alpha_1)\geq 0$, and on contrary, let us also assume $(\beta_2-\beta_1)\geq 0$. This leads to $\alpha_2 \geq \alpha_1$ and $\beta_2 \geq \beta_1$, thereby implying that $\mathcal{P}_1(\alpha_1,\beta_1)$ is isomorphic to some induced subgraph of $\mathcal{P}_2(\alpha_2,\beta_2)$, or both are isomorphic. Hence, utilizing Lemma \ref{ind}, we deduce $\rho(\mathcal{P}_1(\alpha_1,\beta_1))< \rho(\mathcal{P}_2(\alpha_2,\beta_2))$, which is a contradiction. Therefore, two non-isomorphic pineapple graphs $\mathcal{P}_1(\alpha_1,\beta_1)$ and $\mathcal{P}_2(\alpha_2,\beta_2)$ have the same spectral radius then $(\alpha_2-\alpha_1)(\beta_2-\beta_1)<0$.
\end{proof}

\begin{lemma}\label{lemma:2roots}
Let $P(\alpha_1,\beta_1,x)$ and $P(\alpha_2,\beta_2,x)$ be two
polynomials defined by the equation (\ref{char}). Then they have two
common roots if and only if $\beta_i=
\frac{k(\alpha_j-1)(\alpha_j-2)}{(k-1)}~(i\neq j ~~\text{and}
~~i,j=1,2)$ and   $\alpha_1+\alpha_2-2=k$, and the common roots are
given by
$$ \frac{(k-1)\pm \sqrt{(k-1)^2+4(\beta_i-k(\alpha_j-2))}}{2},~ (i,j=1 ~~ \text{or}~~ 2, i\neq j),$$
where $k=\frac{\beta_2-\beta_1}{\alpha_1-\alpha_2}$.
\end{lemma}
\begin{proof}
Let $\rho$ be the common root of the equations
$P(\alpha_1,\beta_1,x)=0$ and $P(\alpha_2,\beta_2,x)=0$, Then, we
have $P(\alpha_1,\beta_1,\rho)=0$ and $P(\alpha_2,\beta_2, \rho)=0$.
After subtracting these equations, we get
 \begin{equation}\label{int4}
 (\alpha_2-\alpha_1)\rho^2+(\alpha_2-\alpha_1+\beta_2-\beta_2)\rho +(\alpha_1\beta_1-\alpha_2\beta_2)+2(\beta_2-\beta_1)=0.
 \end{equation}
 Equation (\ref{int4}) can be rearranged as
  \begin{equation}\label{int5}
 \rho^2-(k-1)\rho -(\beta_1-k(\alpha_2-2))=0 \quad \text{or} \quad  \rho^2-(k-1)\rho -(\beta_2-k(\alpha_1-2))=0,
 \end{equation}
where, $ \quad k= \frac{\beta_2-\beta_2}{\alpha_1-\alpha_2}>0$.
Using division algorithm, we can write $P(\alpha_i,\beta_i,\rho)=0$,
$i=1,2$, as
$$(\rho-\alpha_i+k+1)( \rho^2-(k-1)\rho  -(\beta_i-k(\alpha_j-2))) + A\rho+B_i=0,$$
for $i,j=1,2$ and $i \neq j$, where $A=k(\alpha_1+\alpha_2)-k(k+2)$
and $B=-k(\alpha_i-2)(\alpha_j-2)-(k-1)(\beta_i-k(\alpha_j-2))$.
This implies that  $ A\rho+B_i=0$ if and only if the equations $P(\alpha_i,\beta_i,\rho)=0$, $i=1,2,$ have two common roots, which are given by equation (\ref{int5}).\\
Also, $ A\rho+B_i=0$ further implies that
\begin{equation}\label{eq:A0}
    A=0 \quad \;\Longleftrightarrow\; \quad \alpha_1+\alpha_2=k+2
\end{equation}
and
$$ B=0 \qquad \;\Longleftrightarrow\; -k(\alpha_i-2)(\alpha_j-2)-(k-1)(\beta_i-k(\alpha_j-2))=0.$$
Using equation (\ref{eq:A0}), we obtain
$$ \beta_i= \frac{k(\alpha_j-1)(\alpha_j-2)}{(k-1)}.$$
Clearly, $\alpha_i-k-1=-\alpha_j+1<0$. From equation (\ref{char}),
we know that the product of all the three roots must be negative.
This implies that the quadratic equation (\ref{int5}) is a common
factor of both the equations $P(\alpha_i,\beta_i,\rho)=0$, $i=1,2,$
and have two positive roots. Therefore, the common roots of both the
equations can be  expressed as
$$ \frac{(k-1)\pm \sqrt{(k-1)^2+4(\beta_i-k(\alpha_j-2))}}{2}, $$
$i,j=1 ~~ \text{or}~~ 2$ and $i\neq j$.
\end{proof}

\begin{lemma}\label{lemma:integrs}
Let $k$ and $a$ be two integers with the same parity. Then
$\frac{k(k-a)(k-a-2)}{4(k-1)}$ is an integer if and only if  one of
the following cases hold.
\begin{enumerate}
    {\rm \item [{$(a)$}]  $k \equiv 0 (mod~4) $ \text{and there exists an integer $r$ such that}~~ $a^2-1=(4r+1)(k-1) $,
        \item[$(b)$]  $k \equiv 1 (mod~4)$ \text{and there exists an integer $r$ such that}~~$ a^2-1=(4r+2)(k-1) $,
    \item [{$(c)$}]  $k \equiv 2 (mod~4)$ \text{and there exists an integer $r$ such that}~~ $a^2-1=(4r-1)(k-1)$,
    \item[$(d)$]   $k \equiv 3 (mod~4)$ \text{and there exists an integer $r$ such that}~~ $a^2-1=4r(k-1)$.}
\end{enumerate}
\end{lemma}
\begin{proof}
The quantity $\frac{k(k-a)(k-a-2)}{4(k-1)}$ is an integer if and only if $4(k-1)$ divides $k(k-a)(k-a-2)$, Now, the following cases arise.\\
(a). If $k \equiv 0 (mod~4)$, then
\begin{align}\nonumber
 k(k-a)(k-a-2) &=k\bigg((k-1)^2-2a(k-1)+a^2-1\bigg)\\
               &\equiv a^2-k  (mod~4(k-1))\nonumber\\
               &\equiv (a^2-1)-(k-1)  (mod~4(k-1)).\nonumber
\end{align}
Therefore, there exists an integer $r$ such that
$$(a^2-1)-(k-1) =4r(k-1) \Rightarrow (a^2-1)=(4r+1)(k-1).$$
(b) If $k \equiv 1 (mod~4)$, then
\begin{align}\nonumber
k(k-a)(k-a-2) &=k\bigg((k-1)^2-2a(k-1)+a^2-1\bigg)\nonumber\\
&\equiv a^2-1-2(k-1)  (mod~4(k-1)).\nonumber
\end{align}
Therefore, there exists an integer $r$ such that $a^2-1-2(k-1)=4r(k-1)$, which implies that  $a^2-1=(4r+2)(k-1)$.\\
(c). If $k \equiv 2 (mod~4)$, then
\begin{align}\nonumber
k(k-a)(k-a-2) &=k\bigg((k-1)^2-2a(k-1)+a^2-1\bigg)\nonumber\\
&\equiv a^2-1+(k-1)  (mod~4(k-1)).\nonumber
\end{align}
Therefore, there exists an integer $r$ such that $(a^2-1+(k-1))=4r(k-1)$, which implies that  $a^2-1=(4r-1)(k-1)$.\\
(d). If $k=4m+3$, for some integer $m$, then
\begin{align}\nonumber
k(k-a)(k-a-2) &=k\bigg((k-1)^2-2a(k-1)+a^2-1\bigg)\nonumber\\
&\equiv a^2-1  (mod~4(k-1)).\nonumber
\end{align}
Therefore, there exists an integer $r$ such that $a^2-1=4r(k-1)$.
\end{proof}

Now, we give necessary and sufficient conditions for two
non-isomorphic pineapple graphs to have two common largest
eigenvalues.
\begin{theorem}\label{twocommon}
Let $\mathcal{P}_1(\alpha_1,\beta_1)$ and
$\mathcal{P}_2(\alpha_2,\beta_2)$ be two non-isomorphic pineapple
graphs. Then their two largest eigenvalues are common if and only if
 $\beta_i= \frac{k(\alpha_j-1)(\alpha_j-2)}{(k-1)}$, $i\neq j ~~\text{and} ~~i,j=1,2$,  $\alpha_1+\alpha_2-2=k$, and one of the four conditions of  Lemma \ref{lemma:integrs} holds. The common eigenvalues are given by the expression
 $$ \frac{(k-1)\pm \sqrt{(k-1)^2+4(\beta_i-k(\alpha_j-2))}}{2}, $$
$i,j=1 ~~ \text{or}~~ 2$, $i\neq j$, where
$k=\frac{\beta_2-\beta_1}{\alpha_1-\alpha_2}$ and
$a=\alpha_1-\alpha_2$.
\end{theorem}
\begin{proof}
 Assume that $\alpha_1-\alpha_2=a$. Then, we obtain the following expressions
 $$\beta_1= \frac{k(k-a)(k-a-2)}{4(k-1)}\quad \text{and} \quad  \beta_2= ka-a^2+ \frac{k(k-a)(k-a-2)}{4(k-1)}.
$$
Clearly, $a=\alpha_1-\alpha_2$ and $k=\alpha_1+\alpha_2-2$ have the
same parity. So, Lemma \ref{lemma:integrs} confirms that $\beta_i=
\frac{k(\alpha_j-1)(\alpha_j-2)}{(k-1)}$ are integers if one of the
four conditions in Lemma \ref{lemma:integrs} holds. The rest part
the proof follows by Lemma \ref{lemma:2roots}.
\end{proof}

Now, we have the following observation.
\begin{corollary}\label{cor1}
    Corresponding to two integers $a$ and $k$ of the same parity satisfying one of the four conditions in Lemma  \ref{lemma:integrs}, there exist two pineapple graphs $\mathcal{P}_1(\alpha_1,\beta_1)$ and $\mathcal{P}_2(\alpha_2,\beta_2)$ which have two common largest eigenvalues, where $\alpha_1= \frac{k+a+2}{2} $, $\alpha_2= \frac{k-a+2}{2}$ and  $\beta_i= \frac{k(\alpha_j-1)(\alpha_j-2)}{(k-1)}$, $i\neq j ~~\text{and} ~~i,j=1,2$.
\end{corollary}
\begin{proof}
   Let $a$ and $k$ be two integers with the same parity satisfying one of the four conditions in Lemma \ref{lemma:integrs}. Then, using Lemma \ref{lemma:integrs}, $\alpha_1= \frac{k+a+2}{2} $, $\alpha_2= \frac{k-a+2}{2}$ and  $\beta_i= \frac{k(\alpha_j-1)(\alpha_j-2)}{(k-1)}$ $(i\neq j ~~\text{and} ~~i,j=1,2$ are integers.

Now, we have $\alpha_1-\alpha_2=\frac{k+a+2}{2}-\frac{k-a+2}{2}=a$
and $\alpha_1+\alpha_2= \frac{k+a+2}{2}+\frac{k-a+2}{2}=k+2$. Thus,
according to Theorem \ref{twocommon}, it follows that
$\mathcal{P}_1(\alpha_1,\beta_1)$ and
$\mathcal{P}_2(\alpha_2,\beta_2)$ have two common largest
eigenvalues.
\end{proof}

We illustrate this with the help of the following example.
\begin{example} \end{example}
First choose $a$, say $a=8$. Then we find $a^2-1$  and factorize it
to choose the value of $k$ keeping Lemma \ref{lemma:integrs} in
mind. In our case, we have $a^2-1=63=3\times 21$. Now choose $k=22$.
Then, the pair $(a,k)=(8,22)$ satisfies the third condition of Lemma
\ref{lemma:integrs}. Now, $\alpha_1= \frac{k+a+2}{2}=
\frac{22+8+2}{2}=16 $, $\alpha_2=
\frac{k-a+2}{2}=\frac{22-8+2}{2}=8$ and  $\beta_1=
\frac{k(\alpha_2-1)(\alpha_2-2)}{(k-1)} = \frac{22\times 7\times
6}{(21)}=44$ and $\beta_2= \frac{k(\alpha_1-1)(\alpha_1-2)}{(k-1)} =
\frac{22\times 15\times 14}{(21)}=220$. So, $\mathcal{P}_1(16,44)$
and $\mathcal{P}_2(8,220)$ have two common largest eigenvalues and
the common eigenvalues are 5.783 and 15.217.

\begin{lemma}\label{lemma:root1}
Let $P(\alpha_1,\beta_1,x)$ and $P(\alpha_2,\beta_2,x)$ denote the two polynomials as specified previously. They share precisely one common root $\rho$ if and only if the following conditions hold
\[
\begin{aligned}
    \beta_i = \rho(\rho+1)+\frac{r(\alpha_j-1)}{s}-\frac{r(\rho+1)}{s},  (\alpha_1-\rho-2)(\alpha_2-\rho-2)= \frac{s\rho(\rho+1)}{r},  \text{and} \quad A\neq 0,
\end{aligned}
\]
where $\rho$ is a positive integer, $k=\frac{\beta_2-\beta_1}{\alpha_1-\alpha_2}=\frac{r}{s}$, $r$ and $s$ are coprime, and $A=k(\alpha_1+\alpha_2)-k(k+2)$.

\end{lemma}

\begin{proof}
Proceeding as in Lemma~\ref{lemma:2roots}, if $\rho$ is a common
root of the equations $P(\alpha_i,\beta_i,\rho)=0$ $(i=1,2)$, then
we have $ A\rho+B_i=0$, $i=1,2$, where,
$A=k(\alpha_1+\alpha_2)-k(k+2)$,
$B=-k(\alpha_i-2)(\alpha_j-2)-(k-1)(\beta_i-k(\alpha_j-2))$ and
$k=\frac{\beta_2-\beta_1}{\alpha_1-\alpha_2}$.

 Clearly, $A=0$ gives two common roots. Therefore, for the existence of exactly one common root, we should have $A\neq 0$. This implies that $\rho=-\frac{B}{A}$, which is a rational number. But according to the rational root theorem, a rational root of a monic polynomial is an integer. So, we have $\rho=-\frac{B}{A}$, which gives
\begin{equation}\label{5}
 \frac{k(\alpha_1-2)(\alpha_2-2)+(k-1)(\beta_i-k(\alpha_j-2))}{k(\alpha_1+\alpha_2)-k(k+2)}=\rho~~ \text{(integer)}.
\end{equation}
Let $k=\frac{r}{s}$ such that $\gcd(r,s)=1$. Using the fact that
$\rho$ is also a root of the quadratic equation (\ref{int5}), we
have
\begin{align}
&&\frac{(k-1)\pm\sqrt{(k-1)^2+4(\beta_i-k(\alpha_j-2))}}{2}&=\rho,\nonumber\\
 &or& \quad \frac{(r-s)\pm\sqrt{(r-s)^2+4s(s\beta_i-r(\alpha_j-2))}}{2s}&=\rho,\nonumber\\
 &or& \qquad
 \pm\sqrt{(r-s)^2+4s(s\beta_i-r(\alpha_j-2))}&=2s\rho-(r-s).\nonumber
\end{align}
After squaring and simplification, we obtain
\begin{align}\label{6}
&&(s\beta_i-r(\alpha_j-2))&=\rho(s\rho-r+s)\\
 &or& \qquad \beta_i &= \rho(\rho+1)+\frac{r(\alpha_j-1)}{s}-\frac{r(\rho+1)}{s}, \nonumber
\end{align}
$i,j=1,2$ and $i\neq j$. Using equations (\ref{5}) and (\ref{6}), we
obtain the following expression
$$rs(\alpha_1-2)(\alpha_2-2)+(r-s)(\rho(s\rho-r+s))=rs(\alpha_1+\alpha_2)-r(r+2s)\rho. $$
 Rearranging the terms, we obtain
\begin{equation}\label{alpha}
    (\alpha_1-\rho-2)(\alpha_2-\rho-2)= \frac{s\rho(\rho+1)}{r}.
\end{equation}
\end{proof}

Again, $\alpha_i$ and $\beta_i$ obtained in the previous lemma are
not necessarily integers. Since $\alpha_i$ and $\beta_i$ must be
integers as these are the order of the graphs, the next theorem
deals with the restrictions on the parameters necessary to make
$\alpha_i$ and $\beta_i$ integers.

\begin{theorem}\label{onecommon}

Let $\mathcal{P}_1(\alpha_1,\beta_1)$ and $\mathcal{P}_2(\alpha_2,\beta_2)$ represent two non-isomorphic pineapple graphs. They share exactly one common eigenvalue if and only if the following conditions hold
\[
\begin{aligned}
    \beta_i &= \rho(\rho+1)+\frac{r(\alpha_j-1)}{s}-\frac{r(\rho+1)}{s}, && (\alpha_1-\rho-2)(\alpha_2-\rho-2)= \frac{s\rho(\rho+1)}{r},
    & && \text{and} \quad A\neq 0,
\end{aligned}
\]
where $r$ and $s$ are coprime such that $k=\frac{r}{s}$, $r$ divides $\rho(\rho+1)$, and $s$ divides $\rho+1$, $\alpha_i-1$, for $i=1,2$, and $A=k(\alpha_1+\alpha_2)-k(k+2)$.

\end{theorem}
\begin{proof}
According to Lemma \ref{lemma:root1}, it is sufficient to prove that
$\alpha_i$ and $\beta_i$ are integers. If $\rho$ is the common
eigenvalue, then proceeding as in Lemma \ref{lemma:root1}, we have
$$rs(\alpha_1-2)(\alpha_2-2)+(r-s)(s\beta_i-r(\alpha_j-2))=rs(\alpha_1+\alpha_2)-r(r+2s)\rho. $$
This implies that $ -r^2(\alpha_j-1) \equiv 0 (mod~s)$, which
further implies that  $s| (\alpha_j-1)$, $j=1,2.$ These conditions
along with the expression
$$\beta_i = \rho(\rho+1)+\frac{r(\alpha_j-1)}{s}-\frac{r(\rho+1)}{s}$$
implies that $\beta_i$ are integers if and only if $s$ divides
$\rho+1$. Now, the expression
 $$(\alpha_1-\rho-2)(\alpha_2-\rho-2)= \frac{s\rho(\rho+1)}{r} $$
implies that $\alpha_i$ are integers if and only if
$r|\rho(\rho+1)$,  as $r$ and $s$ are coprime.
\end{proof}

Equation (\ref{alpha}) has as many solutions as the number of
different ways in which $\frac{s\rho(\rho+1)}{r}$ can be factored in
to two factors. But not all the solutions give common spectral
radius. In the next theorem, we show that only those solutions of
equation (\ref{alpha}) can give two graphs
$\mathcal{P}_1(\alpha_1,\beta_1)$ and
$\mathcal{P}_2(\alpha_2,\beta_2)$ having common spectral radius,
which correspond to negative factors of $\frac{s\rho(\rho+1)}{r}$.
In other words, we show that $\alpha_i-\rho-2<0$ or $\alpha_i <
\rho+2$,  for  $i=1,2$.

 \begin{theorem}\label{commonspec}
Let $\mathcal{P}_1(\alpha_1,\beta_1)$ and $\mathcal{P}_2(\alpha_2,\beta_2)$ be two non-isomorphic pineapple graphs. They have exactly one largest common eigenvalue $\rho$ if and only if the following conditions hold
 $$(\alpha_1-\rho-2)(\alpha_2-\rho-2)= \frac{s\rho(\rho+1)}{r} ~~  \beta_i = \rho(\rho+1)+\frac{r(\alpha_j-1)}{s}-\frac{r(\rho+1)}{s} ~~ \text{and}~~ A\neq 0 , $$
where $r$ and $s$ are coprime such that $k=\frac{r}{s}$ and $r$
divides $\rho(\rho+1)$ and $s$ divides $\rho+1$, $\alpha_i-1$ and
$\alpha_i < \rho+2$  $(i=1,2)$, $A=k(\alpha_1+\alpha_2)-k(k+2)$.
 \end{theorem}
 \begin{proof}
Let $\mathcal{P}_1(\alpha_1,\beta_1)$ and
$\mathcal{P}_2(\alpha_2,\beta_2)$ be two non-isomorphic pineapple
graphs. Using Theorem \ref{onecommon}, it is sufficient to prove
that one common eigenvalue of both the graphs is largest if and only
if $\alpha_i < \rho+2$  $(i=1,2)$. Let $\rho$ be the common
eigenvalue of the graphs $\mathcal{P}_1(\alpha_1,\beta_1)$ and
$\mathcal{P}_2(\alpha_2,\beta_2)$. Now, transfer the polynomial in
negative direction by distance $\rho$. Substitute $(y+\rho)$ in
place of $x$ in $P(\alpha_i,\beta_i,x)=0$, we obtain
\begin{equation}\label{eq:shift}
    y(y^2+(3\rho-\alpha_i+2)y+(3\rho^2-2\rho(\alpha_i-2)-\alpha_i-\beta_i+1))=0.
\end{equation}

So, $\rho$ is the largest eigenvalue of the graphs if and only if
the product of the two non-zero roots of equation $(\ref{eq:shift})$
is positive. Thus, we have
     \begin{equation}\label{8}
     3\rho^2-2\rho(\alpha_i-2)-\alpha_i-\beta_i+1>0.
     \end{equation}
     Now, let
     $$(\alpha_1-\rho-2)(\alpha_2-\rho-2)= \frac{s\rho(\rho+1)}{r}=m.n.$$ Then, without loss of generality, we have
     $$\alpha_1=\rho+2+m  \quad \text{and} \quad  \alpha_2= \rho+2+n.$$
     From equation (\ref{8}) for $i=1$, the graphs have the common spectral radius if and only if
     \begin{align}
     3\rho^2-2\rho(\alpha_1-2)-\alpha_1-\beta_1+1&=
                         3\rho^2-2\rho(\rho+2+m)-\rho-2-m-\rho(\rho+1)\nonumber\\
          &+\frac{r}{s}[\rho+1-\rho+2-m+1]+1,\nonumber\\
     or \quad 3\rho^2-2\rho(\alpha_1-2)-\alpha_1-\beta_1+1&= -(m-4)(2\rho+\frac{r}{s})-(\rho-2+m)+1>0\nonumber
     \end{align}
 if and only if $m<0$, which implies that $n<0$, Hence, we have
 $ (\alpha_i-\rho-2)<0 $, which implies that $ \alpha_i< \rho+2 \quad (i=1,2) $.
 \end{proof}

 Now, we have the following observation.
 \begin{corollary}
     If $\rho>2$ is any positive integer, then $\rho$ is the spectral radius of at least one pineapple graph.
 \end{corollary}
\begin{proof}
     Let $\rho>2$ be a positive integer and $s=1$, $r=\frac{\rho(\rho+1)}{2}$. Then the equation
    $$(\alpha_1-\rho-2)(\alpha_2-\rho-2)=\frac{s\rho(\rho+1)}{r}=2$$ has a solution
    $\alpha_1-\rho-2=-1$, that is, $\alpha_1=\rho+1$ and $\alpha_1-\rho-2=-2$, that is, $\alpha_2=\rho$ and the corresponding $\beta_i's$ are $$\beta_1=\rho(\rho+1)+\frac{r(\alpha_2-1)}{s}-\frac{r(\rho+1)}{s}=0, ~~ \beta_2=\rho(\rho+1)+\frac{r(\alpha_1-1)}{s}-\frac{r(\rho+1)}{s}=\frac{\rho(\rho+1)}{2}.$$
    Therefore, $\mathcal{P}(\rho+1,0)$ and $\mathcal{P}(\rho,\frac{\rho(\rho+1)}{2})$ have the same spectral radius $\rho$. Clearly, $\mathcal{P}(\rho+1,0)$ is isomorphic to the complete graph $K_{\rho+1}$ and $\mathcal{P}(\rho,\frac{\rho(\rho+1)}{2})$ is the pineapple graph.
\end{proof}

\begin{example} \end{example}
To create an example satisfying all the conditions of Theorem
\ref{commonspec}, choose a positive integer $\rho=11$ (say). Now
choose $s$ and $r$ such that $s|(11+1)$ and $r|11(11+1)$. So, we
choose $s=2$ and $r=11$. Now, we need to find $\alpha_i$ satisfying
the condition
$(\alpha_1-13)(\alpha_2-13)=\frac{2\times11\times12}{11}=24.$ We
find factors of 24, which are multiples of $s$. Such pair of factors
are $(2,12)$, $(-2,-12)$, $(4,6)$ and $(-4,-6)$. First we choose
$(\alpha_1-13)(\alpha_2-13)=4\times 6$, and we get $\alpha_1=
13+4=17$, $\alpha_2=13+6=19$, $\beta_1=165$ and $\beta_2=154$. We
can check that $\mathcal{P}(17,165)$ and $\mathcal{P}(19,154)$ have
a common eigenvalue 11. But it is not the spectral radius of any
graph. Now, we choose the pair of factors $(-4,-6)$. Then, we obtain
$\alpha_1=7$ $\alpha_2=9$, $\beta_1=110$ and $\beta_2=99$. In this
case $11>
\frac{(\alpha_1-2)+\sqrt{(\alpha_1-2)^2+3(\alpha_1+\beta_1-1)}}{3}
\approx
 8.74$ and $11> \frac{(\alpha_2-2)+\sqrt{(\alpha_2-2)^2+3(\alpha_2+\beta_2-1)}}{3} \approx
 8.10$.   Hence, we obtain the common spectral radius 11 for $\mathcal{P}(7,110)$ and $\mathcal{P}(9,99)$.\\

now, we have the following result.

\begin{theorem}
     If $\beta\geq 8$, then $c(\mathcal{P}(3,\beta))=2\beta+1$.
 \end{theorem}\label{count}
\begin{proof}
     Any induced subgraph of $\mathcal{P}(3,\beta)$ is a pineapple graph with clique number 3 or a star graph $S_{1,n} (n\leq  \beta+1) $ or $K_3$. The star graph $S_{1,n+1}$ is isomorphic to $\mathcal{P}(2,n)$ and $K_3$ is isomorphic to $\mathcal{P}(3,0)$. In other words, we can say that the set of the non empty induced subgraphs of $\mathcal{P}(3,\beta)$ is given by
     $\mathcal{S}(\mathcal{P}(3,\beta))= \{ \mathcal{P}(i,j)| (2\leq i\leq 3, 0\leq j\leq \beta )\}. $
     Let $\mathcal{P}(\alpha_1,\beta_1)$, $\mathcal{P}(\alpha_2,\beta_2)$ $\in \mathcal{S}(\mathcal{P}(3,\beta)) $ have the same spectral radius. Then, without loss of generality, we have $ \alpha_2-\alpha_1=1$. We leave the possibility $ \alpha_2-\alpha_1=0$, since in that case one of the graph becomes the induced subgraph of the other (see Lemma \ref{ind} ).\\
{\it Case 1.} Suppose that both the graphs have two common eigenvalues. Since $\alpha_2=3$ and $\alpha_1=2$, so we have $k=2+3-2=3$. According to Corollary \ref{cor1}, we find only one pair of two graphs having two common eigenvalues, which are $\mathcal{P}(2,3)$ and $\mathcal{P}(3,0)$.\\
{\it Case 2.} Suppose that both the graphs have only one eigenvalue
in common. Then $\alpha_1$ and $\alpha_2$ must satisfy
     $$ (\alpha_1-\rho-2)(\alpha_2-\rho-2)= \frac{s\rho(\rho+1)}{r},$$
     which implies that $$(2-\rho-2)(3-\rho-2)= \frac{s\rho(\rho+1)}{r}, $$
     which gives $\rho= 0, \frac{r+s}{r-s}$. So, now we investigate the existence of pineapple graphs having radius $\frac{r+s}{r-s}$.

   Given that $s|\alpha_i \Rightarrow s=1$, it follows that $\rho= (r+1)/(r-1)$. As $\rho$ is an integer, $(r-1)|(r+1)$, which is possible only when $r=2$ or $3$.

   If $r=2$, then $\rho = 3$. Accordingly, Theorem \ref{commonspec} provides one pair of pineapple graphs having radius $\rho=3$, namely $\mathcal{P}(2,8)$ and $\mathcal{P}(3,6)$.

   If $r=3$, then $\rho=2$, which implies $A=0$. Consequently, no such pair exists in this case. This, combined with Lemma \ref{lem23}, completes the proof.
\end{proof}

\begin{theorem}
\begin{enumerate}
   {\rm \item [{(a)}] For fixed $\alpha \geq 2$, we have
    $ \lim\limits_{\beta \to\infty} r(\mathcal{P}(\alpha,\beta))=1. $
   \item[(b)]  For fixed $\beta \geq 0$, we have
    $ \lim\limits_{\alpha \to\infty} r(\mathcal{P}(\alpha,\beta))=1 $}
\end{enumerate}
\end{theorem}
\begin{proof}
{\bf (a)} For fixed $\alpha\geq 2$, let $\mathcal{S}(\mathcal{P}(\alpha,\beta))$ be the set of the connected induced subgraphs of $\mathcal{P}(\alpha,\beta)$. If any two graphs in $\mathcal{S}(\mathcal{P}(\alpha,\beta))$ have the common spectral radius, then either they have two common eigenvalues or they have only one common eigenvalue.\\
{\it Case 1.} When the graphs $\mathcal{P}(\alpha_1,\beta_1)$ and
$\mathcal{P}(\alpha_2,\beta_2)$ have two common eigenvalues, then
 $\beta_i= \frac{k(\alpha_j-1)(\alpha_j-2)}{(k-1)}$ is fixed for fixed $\alpha_1$ and $\alpha_2$.
This implies that there are finite number of possible pairs of graphs, $\mathcal{P}(\alpha_1,\beta_1)$ and $\mathcal{P}(\alpha_2,\beta_2)$, having two common eigenvalues.\\
{\it Case 2.} The graphs $\mathcal{P}(\alpha_1,\beta_1)$ and
$\mathcal{P}(\alpha_2,\beta_2)$ have exactly one common eigenvalue,
that is, the spectral radius. In this case, they have the common
spectral radius $\rho$
 satisfying the equation
 \begin{equation}\label{1}
    (\alpha_1-\rho-2)(\alpha_2-\rho-2)= \frac{s\rho(\rho+1)}{r}.
 \end{equation}
If $\frac{s}{r}\geq 1$, then one of the negative factors of  $\frac{s\rho(\rho+1)}{r}$ will be lesser than $-(\rho+1)$. Without loss of generality, we have $\alpha_1-\rho-2\leq-\rho-1$. This gives $\alpha_1\leq 1$, which is a contradiction. \\
If $\frac{s}{r} < 1$, then equation (\ref{1}) can be written as
 $$ \left(1-\frac{s}{r}\right)\rho^2-\rho\left(\alpha_1+\alpha_2-4-\frac{s}{r}\right)+\alpha_1\alpha_2-2(\alpha_1+\alpha_2)+4=0.$$

 If $\rho_1$ and $\rho_2$ are the two roots of this quadratic equation, then
 \begin{align*}
 \rho_1+\rho_2&=\frac{\alpha_1+\alpha_2+\dfrac{s}{r}-4}{1-\dfrac{s}{r}}= (\alpha_1+\alpha_2-4)+\dfrac{s(\alpha_1+\alpha_2-3)}{(r-s)}\\
 &\leq (\alpha_1+\alpha_2-3)(s+1)-1\leq \alpha(2\alpha-3)-1.
 \end{align*}
This means that, for fixed $\alpha$, the spectral radius $\rho$ is
bounded above, so are $\beta_i$, $i=1,2$ as
$$\beta_i = \rho(\rho+1)+\frac{r(\alpha_j-1)}{s}-\frac{r(\rho+1)}{s}= \rho(\rho+1)+\frac{r(\alpha_j-\rho-2)}{s}\leq \rho(\rho+1).$$
 The last inequality is due to the fact that $\alpha_j-\rho-2<0$.
This implies that there are finite number of possible pairs of
graphs $\mathcal{P}(\alpha_1,\beta_1)$ and
$\mathcal{P}(\alpha_2,\beta_2)$ having the common spectral radius.
Thus, there exists $\beta^*(\alpha)$ depending on $\alpha$ such that
the quantity
$b(\mathcal{P}(\alpha,\beta))-c(\mathcal{P}(\alpha,\beta))$ is
constant for all $\beta \geq \beta^*(\alpha)$. This gives
 $$ \lim\limits_{\beta \to\infty} r(\mathcal{P}(\alpha,\beta))= \lim\limits_{\beta \to\infty} \bigg(1+\frac{b(\mathcal{P}(\alpha,\beta))-c(\mathcal{P}(\alpha,\beta))}{c(\mathcal{P}(\alpha,\beta))}\bigg) =1$$
  as $c(\mathcal{P}(\alpha,\beta))$ increases unboundedly when $\beta \to \infty$.\\
{\bf (b)} Similarly, for fixed $\beta\geq 0$, we have the following cases.\\
{\it Case 1.} When the graphs $\mathcal{P}(\alpha_1,\beta_1)$ and
$\mathcal{P}(\alpha_2,\beta_2)$ have two common eigenvalues. Then,
we have
  $\alpha_1+\alpha_2-2=k=\frac{\beta_2-\beta_1}{\alpha_1-\alpha_2} \leq \beta$. This gives $\alpha_i \leq \beta+2$, which implies that there are finite number of pairs of graphs $\mathcal{P}(\alpha_1,\beta_1)$ and $\mathcal{P}(\alpha_2,\beta_2)$ having two common eigenvalues.\\
{\it Case 2.}  When the graphs $\mathcal{P}(\alpha_1,\beta_1)$ and
$\mathcal{P}(\alpha_2,\beta_2)$ have exactly one common eigenvalue,
that is, the spectral radius. In this case, we have $r \leq
\frac{\beta_2-\beta_1}{\alpha_1-\alpha_2}s\leq \beta $,
  \begin{align}
  \beta_i &= \rho(\rho+1)+\frac{r(\alpha_j-1)}{s}-\frac{r(\rho+1)}{s}  \geq \rho(\rho+1)-\frac{r(\rho+1)}{s}\nonumber\\
  & \geq \rho(\rho+1)-\beta(\rho+1).\nonumber
  \end{align}
 This implies that
$ \beta  \geq \rho(\rho+1)-\beta(\rho+1)$, which further gives
  $\beta \geq \frac{\rho(\rho+1)}{(\rho+2)}\geq \rho-1.$ This implies that the spectral radius is bounded above, so are $\alpha_i$ as $\alpha_i< \rho-2$. Therefore, there are finite pairs of graphs having the common spectral radius. So, there exists $\alpha^*(\beta)$ depending on $\beta$ such that the quantity $b(\mathcal{P}(\alpha,\beta))-c(\mathcal{P}(\alpha,\beta))$ is constant for all $\alpha \geq \alpha^*(\beta)$. This gives
$$ \lim\limits_{\alpha \to\infty} r(\mathcal{P}(\alpha,\beta))= \lim\limits_{\alpha \to\infty} \bigg(1+\frac{b(\mathcal{P}(\alpha,\beta))-c(\mathcal{P}(\alpha,\beta))}{c(\mathcal{P}(\alpha,\beta))}\bigg) =1$$
as $c(\mathcal{P}(\alpha,\beta))$ increases unboundedly when $\beta \to \infty$.\\
\end{proof}

\section*{Acknowledgements}
The research of Pawan Kumar is supported by CSIR, India as a Senior
Research Fellowship, file No. 09/112(0669)/2020-EMR-I. The research of S. Pirzada is supported by NBHM-DAE research project number NBHM/02011/20/2022.

\section*{Declaration}

The authors report there are no competing interests to declare

\section*{Data availability statement}

The paper does not use or produce any dataset.


\begin{thebibliography}{99}
\bibitem{Alon}N. Alon, B. Bollobas, Graphs with a small number of distinct induced subgraphs, \textit{Discrete Math.} {\bf  75} (1989) 23-30.
\bibitem{Fer}R. Fernandes, J. Judice, V. Trevisan, Complementary eigenvalues of graphs. \textit{Linear Algebra  Appl.} {\bf 527} (2017) 216-231.
\bibitem{Fer2} L. Fernandes, J. Judice, H. Sherali, M. Fukushima, On the computation of all eigenvalues for the eigenvalue complementarity problem, \textit{J. Global Optim.} {\bf 59} (2014) 307-326.
\bibitem{ppmt} S. Merajuddin, Pawan Kumar, S. Pirzada, Vilmar Trevisan, A unified criterion for distinguishing graphs by their spectral radius, \textit{Linear Multilinear Algebra} to appear. https://doi.org/10.1080/03081087.2023.2228458
\bibitem{Pin2020} L. K. Pinheiro, B. S. Souza, V. Trevisan, Determining graphs by the complementary spectrum, \textit{Discussiones Mathematicae Graph Theory} {\bf 40} (2020) 607-620.
\bibitem{See2020} A. Seeger, Repetition of spectral radii among connected induced subgraphs, \textit{Graphs Combinatorics} {\bf 36} (2020) 1131-1144.
\bibitem{Sos2021}A. Seeger, D. Sossa, On cardinality of complementarity spectra of connected graphs, \textit{Linear Algebra Appl.} {\bf 614} (2021) 5-23.
\bibitem{See2018}A. Seeger, Complementarity eigenvalue analysis of connected graphs, \textit{Linear Algebra Appl.} {\bf 543} (2018) 205-225.
\bibitem{intr} A. Seeger, Eigenvalue analysis of equilibrium processes defined by linear complementarity conditions, \textit{ Linear Algebra Appl.} {\bf 292} (1999) 1-14.
\bibitem{see2}A. Seeger, D. Sossa, Extremal problems involving the two largest complementarity eigenvalues of a graph, \textit{Graphs Combinatorics} {\bf 36} (2020) 1-25.
\bibitem{friendship} A. Seeger, D. Sossa, Spectral radii of friendship graphs and their connected induced subgraphs, \textit{Linear Multilinear Algebra} {\bf } {\bf 71(1)} (2023) 63-87.
\bibitem{sha} Y. Shang, Groupies in multitype random graphs.
\textit{SpringerPlus} {\bf 5} (2016) 989.
\bibitem{ineq} Z. Stanic, Inequalities for Graph Eigenvalues, Cambridge University Press, 2015.
\bibitem{spec} D. Stevanovic,  Spectral Radius of Graphs. In: A. Encinas, M. Mitjana (eds) Combinatorial Matrix Theory. Advanced Courses in Mathematics-CRM Barcelona. Birkhäuser, Cham. (2018)
\bibitem{topcu2016spectral}H. Topcu, S. Sorgun, W. H. Haemers, On the spectral characterization of pineapple graphs, \textit{Linear Algebra Appl.} {\bf 507} (2016) 267-273.
\bibitem{topcu2019graphs} H. Topcu, S. Sorgun, W. H. Haemers, The graphs cospectral with the pineapple graph, \textit{Discrete Appl. Math.} {\bf 269} (2019) 52-59.
\bibitem{zhang2009some} X. Zhang, H. Zhang, Some graphs determined by their spectra, \textit{Linear Algebra  Appl.} {\bf 431} (2009) 1443-1454.

\end{thebibliography}
\end{document}